# DISCUSSION OF "BREAKDOWN AND GROUPS" BY P. L. DAVIES AND U. GATHER

By Xuming He

*University of Illinois at Urbana–Champaign*

The notion of breakdown point has influenced the robust statistics literature for over three decades. Professors Davies and Gather make a convincing argument that the common understanding of a high breakdown point is intimately connected to a group structure in the sample space. I would like to applaud the authors for a fine piece of work to formalize the connection. Inspired by their work, I would like to offer my opinions on the nature and the future of the breakdown point.

A high breakdown point is usually considered to be a virtue of a statistical procedure, because such a procedure is less affected by least favorable configurations of data contamination. Thus arises a natural question of how high the breakdown point can be in a given problem. For location equivariant functionals, $1/2$ is a tight upper bound on the breakdown point. In more structured problems and in more general settings, the defintion of a breakdown is less straightforward. Stromberg and Ruppert (1992) considered nonlinear regression. He and Simpson (1992) provided a definition of breakdown for general parameter spaces which might be compact. Instead of citing more work on breakdown in specific settings, I would emphasize that it is the simplicity and intuition that has made the breakdown point a popular measure of global robustness. Intuitively speaking, the breakdown point is the smallest fraction of data contamination that could make an estimator or test statistic totally uninformative or unusable.

The point I would like to make is that it is better to remember the spirit, not the letter, of any definition of the breakdown point. To illustrate this point, let us use the following definition of a finite-sample breakdown point.

Given a sample $X_n$ of size $n$, the breakdown point of $T_n(X_n)$ is

$$\varepsilon_n^* = \min\left\{ m/n : \sup_{X_{n,m}^*} d(T_n(X_{n,m}^*), T_n(X_n)) = \infty \right\},$$

---

where $X_{n,m}^*$ is obtained by replacing $m$ out of $n$ points in $X_n$ by arbitrary values, and $d(a,b) \in (0,\infty)$ is some distance measure between $a$ and $b$. Let us assume here that $d(a,b)$ can take arbitrarily large values. If $T$ is a location estimator, one often takes $d(a,b) = |a-b|$. If $T$ is a scale estimator, one may take $d(a,b) = |\log(a/b)|$. For most location, scale and regression estimators, the breakdown points do not depend very much on the initial sample $X_n$, but this is not always the case. Because the breakdown point is defined at each sample, we can easily modify any estimator $T_n$ so that it will not break down at all according to this definition. For example, we can take $T_n^*(X_n) = \max\{-n, \min\{n, T_n(X_n)\}\}$ for a location estimator, and it will never be unbounded at any contamination.

Such a construction, however, violates the spirit of a high breakdown estimator, albeit it is mathematically legitimate. For a location estimator this problem can be eliminated by imposing location equivariance. In a general setting it is not clear what can be done. I use this example to stress that we should not try to exploit the mathematics of a statistical concept without a clear sense of purpose.

When someone claims to have found an estimator with breakdown point equal to 1, my first reaction tends to be that it might not be an appropriate use of the notion. Understanding and imposing a group equivariance structure on the estimator certainly helps, but it cannot eliminate inappropriate use of breakdown. In some problems (e.g., logistic regression) the group structure that can be identified might be very limited.

The notion of breakdown for a test statistic does not always carry the same implications as for an estimator. Davies and Gather discuss in their treatment of logistic regression whether the parameter value of 0 should be considered as a breakdown. I agree with the authors that the value of 0 plays no special role for an estimator. To study the breakdown of a statistical test, the value 0 often plays a special role. I refer to He, Simpson and Portnoy (1990) for more detail, but simply point out the obvious that one cannot judge the appropriateness of a breakdown definition without further specifics.

Take, for example, the classification tree. It is reasonable to say that a procedure breaks down if the classification rule becomes no better than a random guess. However, it is not obvious at all how to construct a tree with the highest possible breakdown point. Other notions of breakdown are also possible here.

I hope that the breakdown point will remain as a simple and intuitive concept. Maybe it falls into the same category as "outlier," where some degree of vagueness would win over more users. When every statistician starts to talk about his or her own notion of a breakdown point, I think that we have made it.

Having made my main point, I would like to use this opportunity to offer my thoughts on some of the controversial issues surrounding the breakdown point.

1. Is the breakdown point a conservative measure of robustness? Yes, it is by definition. But as long as we know what it is doing, it is not bad to be conservative. 2. Are there good reasons to aim for the highest possible breakdown point? Usually not. I am not voicing objections to research on the highest possible breakdown point procedures; such research can offer insights. In choosing a good statistical procedure, we have to balance breakdown with other measures of quality. 3. Some say that high breakdown estimators are usually locally unstable. There is some truth to this, but again, one has to strike a balance between breakdown and local stability. This statement would be as true or untrue as "efficient estimators usually have low breakdown points." 4. High breakdown point estimators are usually too hard to compute. It is easy to propose a difficult-to-compute high breakdown estimator, but advances in methodological research and in computing power are already making more and more high breakdown procedures practical. Obviously I like the fact that SAS procedures based on high breakdown method are being added.

Finally, what role will the notion of breakdown point play in the future? I am not good at predicting the future, but I hope that it will be in every statistician's mind in evaluating the quality of a statistical procedure. It is in our best interest to keep it as simple and intuitive as possible so that it will be understood and appreciated by every statistician (plus more). In addition to research papers such as this one under discussion, I hope to see educational papers, too, that will be accessible by a broader audience. If I use the NSF jargon, I hope to see both scientific merit and broader impacts. I think that we will get there if we all try.

Department of Statistics
University of Illinois at
 Urbana–Champaign
725 S. Wright Street
Champaign, Illinois 61820
USA
e-mail: x-he@uiuc.edu